\documentclass[reqno,12pt]{article}
\usepackage{hyperref}
\usepackage{amsmath}
\usepackage{amssymb}
\usepackage{mathrsfs}

\newtheorem{thm}{Theorem}[section]
\newtheorem{lemma}[thm]{Lemma}

\newtheorem{corol}[thm]{Corollary}
\newtheorem{propos}[thm]{Proposition}
\newtheorem{rema}{Remark}[section]
\def\bp{\begin{propos}}
\def\ep{\end{propos}}
\def\bt{\begin{thm}}
\def\et{\end{thm}}
\def\bco{\begin{corol}}
\def\eco{\end{corol}}
\def\bl{\begin{lemma}}
\def\el{\end{lemma}}
\def\br{\begin{rema}}
\def\er{\end{rema}}
\def\be{\begin{equation}}
\def\ee{\end{equation}}
\def\ba{\begin{array}}
\def\ea{\end{array}}
\def\bena{\begin{eqnarray}}
\def\eena{\end{eqnarray}}

\def\P{{\mathbb P}}
\def\E{{\mathbb E}}

\def\1{I}

\def\a{{\alpha}}
\def\D{{\Delta}}

\def\var{\hb{Var}}

\def\QED{\hfill$\square$\vskip 3mm}

\def\Dp{\displaystyle}
\def\Df{\Dp\frac}
\def\hb{\hbox}
\def\({\left(}
\def\){\right)}
\def\[{\left[}
\def\]{\right]}

\begin{document}
\title{\LARGE The Degree Sequence of a Scale-Free Random Graph\\
Process with Hard Copying\\[5mm]
 \footnotetext{Key
words and phrases: degree sequence; power law; hard copying; random
graph process}}
\author{ Gao-Rong Ning$^1$\thanks{Supported in part
by the Natural Science Foundation of China},\ \ Xian-Yuan Wu$^{1*}$\
\hb{ and
 }Kai-Yuan Cai$^2$\thanks{Supported in part by the Natural
Science Foundation of China and MicroSoft Research Asia under grant
60633010} }

\date{}

\maketitle {\small \vskip -5mm
\begin{center}
\begin{minipage}{13cm} \noindent\hskip -2mm$^1$School of
Mathematical Sciences, Capital Normal University,
Beijing, 100037, China. Email: \texttt{ninggaorong@163.com}; \texttt{wuxy@mail.cnu.edu.cn}\\[-5mm]

\noindent\hskip -2mm$^2$Department of Automatic Control, Beijing
University of Aeronautics and Astronautics, Beijing, 100083, China.
Email: \texttt{kycai@buaa.edu.cn}
\end{minipage}
\end{center}
\vskip 5mm
\begin{center} \begin{minipage}{15cm}
{\small {\bf Abstract}: In this paper we consider a simple model of
random graph process with {\it hard} copying as follows: At each
time step $t$, with probability $0<\alpha\leq 1$ a new vertex $v_t$
is added and $m$ edges incident with $v_t$ are added in the manner
of {\it preferential attachment}; or with probability $1-\alpha$ an
existing vertex is copied uniformly at random. In this way, while a
vertex with large degree is copied, the number of added edges is its
degree and thus the number of added edges is not upper bounded. We
prove that, in the case of $\alpha$ being large enough, the model
possesses a mean degree sequence as $ d_{k}\sim Ck^{-(1+2\alpha)}$,
where $d_k$ is the limit mean proportion of vertices of degree $k$.
}
\end{minipage}
\end{center}
}

\vskip 5mm
\section{Introduction and the statement of the main result}
\renewcommand{\theequation}{1.\arabic{equation}}
\setcounter{equation}{0}

Real-world networks such as economic companies, biological
oscillators, social networks, and the World Wild Web (internet) {\it
etc.} can be modeled by random complex graphs
\cite{BKEMS,KRRS,LLJ,M,St,WS}. By studying random complex graphs,
various topological properties such as degree-distribution
\cite{BA,BO,BRST,CFV}, diameter \cite{ABJ,ASBS,BR2}, clustering
\cite{BR,N}, stability \cite{B,B2,BR3} and spectral gap \cite{ACL}
of these real-world networks have been presented. One of the most
basic properties of many real-world networks is concerned with the
power law degree distributions. As indicated in \cite{BA}, the
emergence of the power law degree distributions should be a
consequence of two generic mechanisms:
\begin{enumerate}
             \item Evolution: new vertices and edges are added continuously, and
             \item Preferential attachment: new vertices are preferentially attached to vertices that are already well
connected,
\end{enumerate}
The above mechanisms are referred to as BA mechanisms. Besides the
original model proposed in \cite{BA}, many other models with the BA
mechanisms have been introduced and aimed to explain the underlying
causes for the emergence of the power law degree distributions. This
can be observed in `LCD model' \cite{BR2}, the generalization of
`LCD model' due to Buckley and Osthus \cite{BO}, the very general
models defined by Copper and Frieze \cite{CF}, Copper, Frieze and
Vera \cite{CFV} {\it etc.}



Copying is another mechanism that may be observed in real-world
networks. The basic idea of copying comes from the fact that a new
web page is often made by copying an old one. A kind of copying
models was proposed in Kumar {\it et al}. \cite{KRRS} to explain the
emergence of the degree power laws in the web graphs. These models
are parameterized by a {\it copy factor} $\a\in (0,1)$ and a
constant out-degree $d\geq 1$. At each time step, one vertex $u$ is
added and $d$ out-links are generated for $u$ as follows. First, an
existing vertex $p$ is chosen uniformly at random; then with
probability $1-\a$ the $i^{\rm th}$ out-link of $p$ is taken to be
the $i^{\rm th}$ out-link of $u$, and with probability $\a$ a vertex
is chosen from the existing vertices uniformly at random to be the
destination of the $i^{\rm th}$ out-link of $u$. It is proved in
\cite{KRRS} that the above copying models possess a power law degree
sequence as $d_k\sim Ck^{-(2-\a)/(1-\a)}$.


In this paper we will introduce and study a new copying model
created by {\it lazy} copiers. Our copiers are so lazy that the only
thing they want to do is copying. However, the copiers corresponding
to the copying action discussed in \cite{KRRS} should be more clever
and diligent: for the chosen vertex $p$, they have to distinguish
which link be a original out-link of $p$ first and then decide
whether or not to copy it.

Let's consider the following random process $G_{t}$, $t = 2, 3,
\cdots$. Assume that graph $G_{t} = (V_{t},E_{t})$ and $t=|V_t|$, $
e_{t}=|E_t|$ (In order to simplify the statement and the proof of
our main result, technically, we start our process at time step 2).
\vskip 3mm Time-Step 2: To begin the process, we start with $G_2$
consisting of vertices $v_{1}$, $v_2$ and $2m$ multi-edges between
them. \vskip 3mm

Time-Steps $ t\geq 3$:
\begin{itemize}
\item With probability $\a>0$ we add a new vertex $ v_{t}
$ to $G_{t-1}$ and then add $m$ random edges incident with $v_{t}$.
The $m$ random neighbors $ w_{1}, w_{2}, \ldots, w_{m}$ are chosen
independently and for any $1 \leq i \leq m$, $w \in V_{t-1}$,
\be\label{1.1}\P(w_{i} = w) = \Dp\frac{d_{w}(t-1)}{2e_{t-1}},\ee
where $d_w(t-1)$ denotes the degree of vertex $w$ in $G_{t-1}$. Thus
neighbors are chosen by {\it preferential attachment.}
\item With probability $1-\a $ we generate vertex $ v_{t} $
by copying a existing vertex $v_i$, $1\leq i\leq t-1$ from $V_{t-1}$
uniformly at random. Note that in this case, all neighbors of $v_t$
are those of the copied vertex $v_{i}$.
\end{itemize}
As defined above, our copying is executed in a direct and simple
way, which is referred to as {\it hard} copying here. With hard
copying, $e_t$ may increase nonlinearly, this makes bounding $e_t$ a
rather hard problem.

Now, Let $D_k(t)$ be the number of vertices with degree $k\geq 0$ in
$G_t$ and let $\overline{D}_k(t)$ be the expectation of $D_k(t)$.
The main result of this paper follow as:

\bt\label{th1} Assume that $2m(1-\a)<\a$. Then, for all $k\geq 0$,
the limit
$d_k=\Dp\lim_{t\rightarrow\infty}\frac{\overline{D}_k(t)}{t}$ exists
and satisfies $$ d_k=0,\ 0\leq k<m;\ d_m=\Df{2\a}{m+2\a};\
d_k=\prod_{i=m+1}^k\(1+\Df{1+2\a}{i+2\a}\)d_m,\ \forall\ k>m.$$
Obviously, $d_k\sim Ck^{-(1+2\a)}$ for some constant $C$. \et

We follow the basic procedures in \cite{CF} and \cite{CFV} to prove
our main theorem. The rest of the paper is organized as follows. In
Section 2, we bound the maximum degree and then bound $e_t$, the
number of edges in $G_t$. In Section 3, using the estimates given in
Section 2, we establish the recurrence for $\overline{D}_k(t)$.
Finally, in section 4, we derive the approximation of
$\overline{D}_k(t)$ by a recurrence with respect to $k$ and then
solve the recurrence in $k$ to finish the proof of
Theorem~\ref{th1}.

Here we note that although this paper focuses on the power law
degree distributions, other degree distributions including the
exponential degree distributions of random graph process have also
been observed \cite{ASBS,BKEMS,LLJ,WS}. Furthermore, phase
transition may emerge in the degree distributions of random graph
processes \cite{WDLC,WDLC2}. The phase transition problem of the
copying model proposed in this paper is left to future
investigation.

\vskip 5mm
\section{Bounding the degree and the number of edges }
\renewcommand{\theequation}{2.\arabic{equation}}
\setcounter{equation}{0}

In this section, we first bound the maximum degree in $G_t$ and then
bound $e_t$. Actually, we will give four kinds of estimates to
$e_t$, as will be seen in section 3, the four estimates are all
necessary for establishing the recurrence of $\overline{D}_k(t)$.

For $t\geq 2$, let $V^o_t$ be set of {\it original } vertices in
$V_t$, namely
$$V^o_t:=\{v\in V_t: v=v_1,\ v_2\hb{ or }v\hb{ is added as a new vertex at
some time step }3\leq s\leq t\}.$$ For any times $s$ and $t$ with
$3\leq s\leq t$, if $v_s\in V^o_t$, then, \be\label{2.3}
d_{v_s}(s)=\frac 12d_{v_1}(2)=\frac12d_{v_2}(2)=m.\ee

We say an event happens {\it quite surely} (qs) if the probability
of the complimentary set of the event is $O(t^{-K})$ for any $K>0$.

We bound the degree in $G_t$ from top as follows

\bl\label{l1} Assume that $2m(1-\a)<1$ and $v_s\in V^o_t$. Then
\be\label{1}d_{v_{s}}(t) \leq \({t}/{s}\)^{{\a}/{2}+m(1-\a)}(\log
t)^{3}\ \ \  qs. \ee \el

{\it Proof}: Let $Y$ be the $\{0,1\}$-valued random variable with
$\P(Y=1)=\a=1-\P(Y=0)$. Then using the fact that $e_t\geq mt$, we
have \be\label{2.2} \E(d_{v_{s}}(t+1)\mid G_t)\leq
d_{v_{s}}(t)+YB\(m,\frac{d_{v_{s}}(t)}{2mt}\)+(1-Y)mB\(1,\frac{d_{v_{s}}(t)}{t}\),\ee
where $B(\cdot,\cdot)$ be the general Binomial random variable.

Using the fact (\ref{2.3}) and the relation (\ref{2.2}),
Lemma~\ref{l1} follows from the same argument as used in \cite{CF},
\cite{CFV} and \cite{WDLC}.\QED

For any $v\in V_t$, if $v$ is copied at time step $s$ from some
vertex $v_r$, $1\leq r\leq s-1$, we call $v$ the {\it daughter}
vertex of $v_r$ and call $v_r$ the {\it mother} vertex of $v$.
Denote by $D(v,G_t)$ the set of all descendants of $v$ in $G_t$. By
the definition of the model, we know that, for any $v_s\in V^o_t$
and $v\in D(v_s,G_t)$, $d_v(t)$ is same distributed as $d_{v_s}(t)$.
Now, denote by $\Delta_{t}$ the maximum degree in $G_{t}$, then, by
Lemma~\ref{l1} and the above analysis, we have
\be\label{2.5}\Delta_{t}\leq t^{{\a}/{2}+m(1-\a)}(\log t)^{3},\ \
qs.\ee

For any $v_s\in V^o_t$, let $f_{v_{s}}(t)=|D(v_s,G_t)|$ be the
number of all descendants of $v_s$, then, we have

\bl\label{l2}For any $s\geq 1$, if $v_s$ is a original vertex, i.e.,
for some $t\geq 2$, $v_s\in V^o_t$, then
\be\label{2.6}f_{v_s}(t)\leq \(t/s\)^{1-\a}\(\log t\)^3,\ \ \
qs.\ee\el

{\it Proof}: Let $Y$ be the random variable used in the proof of
Lemma~\ref{l1}, then,

\be\label{2.7} \E(f_{v_{s}}(t+1)\mid
G_t)=f_{v_{s}}(t)+(1-Y)B\(1,\frac{f_{v_{s}}(t)}{t}\).\ee The Lemma
follows from the relation (\ref{2.7}) and the same argument as used
in Lemma~\ref{l1}.\QED

Now we begin to bound $e_t$, the number of edges in $G_t$. Let
$a_{t}$ be the number of edges added at time step $t+1$, i.e.,
$e_{t+1}=a_{t}+e_{t}$. By the definition of the model, we have
$a_t\leq \max\{\Delta_t,m\}=\D_t$, $\forall\ t\geq 2$; on the other
hand, noticing that the number of multi-edges between any given
vertices pair is fewer than $2m$, we have
$$ \Delta_2=2m,\ \ \Delta_{t+1}\leq\Delta_t+2m,\ \forall\ t\geq
2.$$ This gives the following determined upper bound on $e_t$,
\be\label{2.22}e_t=2m+\Dp\sum_{s=2}^{t-1}a_s\leq
2m+\Dp\sum_{s=2}^{t-1}2m(s-1)=O(t^2).\ee

For random upper bounds on $e_t$, firstly, we prove a crude one as
\be\label{2.21}e_t\leq O\(t(\log t)^6\),\ \ qs.\ee Indeed, we have
$$2e_t=\sum_{s=1}^t d_{v_s}(t)=\sum_{v_s\in V^o_t}\sum_{v\in
D(v_s,G_t)}d_v(t).$$By Lemma~\ref{l1} and Lemma~\ref{l2}, \bena
\sum_{v_s\in V^o_t}\sum_{v\in D(v_s,G_t)}d_v(t)\leq \sum_{s=1}^t
\[\(t/s\)^{\a/2+(m+1)(1-\a)}(\log t)^6\]=O\(t(\log t)^6\),\ \
qs.\nonumber\eena Note that for the last equality we have used the
condition $2m(1-\a)<\a$, which is given in the statement of
Theorem~\ref{th1}.

Secondly, we try to give an estimate to $\E(e_t)$, the expectation
of the number of edges in $G_{t}$. By the definition of the model,
we have \be\label{2.11}\E(e_{t+1}|G_{t})=e_{t}+\a
m+(1-\a)\frac{2e_{t}}{t},\ee so
\be\label{2.8}\E(e_{t+1})=\E(e_{t})\(1+\frac{2(1-\a)}{t}\)+\a m.\ee
Let
$$\eta_{t}:=e_{t}-\mu t,$$ where $\mu=\Dp\frac{\a m}{1-2(1-\a)}$.
Then, (\ref{2.8}) implies that
$$\E(\eta_{t+1})=\E(\eta_{t})\(1+\frac{2(1-\a)}{t}\).$$
Thus, $\E(\eta_{t})=O(t^{2(1-\a)})$ and we have
\be\label{2.9}\E(e_{t})=\mu t+O(t^{2(1-\a)}).\ee

Finally, we have the following probability estimate on $e_t$ as
\bl\label{l3} Assume that $2m(1-\a)<1$. Take $\varepsilon_0
>0$ such that $1+2\varepsilon_0+2m(1-\a)<2$, then
\be\label{2.10}\P\(|e_{t}-\mu
t|>t^{\frac{1}{2}+\varepsilon_0+m(1-\a)}\)
=O(t^{-\varepsilon_0}).\ee\el

{\it Proof}: To get the estimate (\ref{2.10}), we have to bound
$\var (e_t)$, the variance of $e_t$. First of all, we have
\be\label{2.17} \var(e_{t+1})= \var(a_{t}+e_{t})=
\var(e_{t})+\var(a_{t})+2\(\E(a_{t}e_{t})-\E(a_{t})\E(e_{t})\).\ee

By definition, we have $$ \E(a^2_t\mid G_t)=\a
m^2+(1-\a)\sum_{s=1}^t\Df{d^2_{v_s}(t)}{t}.$$ Then, by
Lemma~\ref{l1} and Lemma~\ref{l2}, \bena \E(a_{t}^{2})&&\hskip-5mm=
\a m^{2}+\Df{(1-\a)}t\E\(\sum_{v_s\in V^o_t}\sum_{v\in
D(v_s,G_t)}d^2_v(t)\)\nonumber\\[3mm]
&&\hskip-5mm\leq \a
m^{2}+\Df{(1-\a)}t\sum_{s=1}^t\[(t/s)^{\a+2m(1-\a)}(\log
t)^6\]\[(t/s)^{1-\a}(\log t)^3\]+O(t^{-10})\nonumber\\[3mm]
&&\hskip-5mm =O\(t^{2m(1-\a)}(\log t)^{9}\).\label{2.12}\eena In
addition, by (\ref{2.11}) and (\ref{2.9}), we have
\be\label{2.13}\E(a_t)=\a m+2(1-\a)\mu+O(t^{2(1-\a)-1}).\ee Thus
\be\label{2.14}\var(a_t)=O\(t^{2m(1-\a)}(\log t)^9\).\ee

For the term $\E(a_te_t)$, using (\ref{2.11}), it is clear that
$$\E(a_{t}e_{t}|G_{t})=e_{t}\E(a_{t}|G_{t})=e_{t}\(m\a+2(1-\a)\frac{e_{t}}{t}\),$$
then
\be\label{2.15}\E(a_{t}e_{t})=m\a\E(e_{t})+\frac{2(1-\a)}{t}\E(e_{t}^{2}).\ee
Using (\ref{2.11}) again, we have \be\label{2.16}\E(a_t)\E(e_t)=m\a
E(e_t)+\frac{2(1-\a)}{t}\E(e_{t})^{2}.\ee

Substituting (\ref{2.14}), (\ref{2.15}) and (\ref{2.16}) into
(\ref{2.17}), we get
\bena\label{2.18}\var(e_{t+1})&&\hskip-5mm=\(1+\Df{4(1-\a)}{t}\)\var(e_t)+O\(t^{2m(1-\a)}(\log
t)^9\)\nonumber\\[3mm]
&&\hskip-5mm=\(1+\Df{4(1-\a)}{t}\)\var(e_t)+O\(t^{2m(1-\a)+\varepsilon_0}\),\label{2.18}\eena
where $\varepsilon_0>0$ is given in the statement of the Lemma. The
recurrence (\ref{2.18}) can be solved directly to get
$$\var(e_t)=\prod_{s=3}^{t-1}\(1+\Df{4(1-\a)}{s}\)\(\var(e_3)+O\(\sum_{s=3}^{t-1}\Df{s^{2m(1-\a)+\varepsilon_0}}{\prod_{j=3}^s\(1+\Dp{4(1-\a)}/{j}\)}\)\)$$
for large $t$, this implies that
\be\label{2.19}\var(e_t)=O\(t^{1+2m(1-\a)+\varepsilon_0}\).\ee The
Lemma follows immediately from (\ref{2.9}), (\ref{2.19}) and the
Chebychev's inequality. \QED \vskip 5mm

\section{Establishing The Recurrence for $\overline{D}_k(t)$}
\renewcommand{\theequation}{3.\arabic{equation}}
\setcounter{equation}{0}

Before we establish the recurrence for $\overline{D}_k(t)$, we have
to bound the multi-edges first. For $t\geq 2$, let $$Z_{t}=\{v\in
V_{t}:\exists\ u\in V_{t}\hb{ s.t. there are multi-edges between } u
\ \hb{ and } v\}$$ and $X_t=|Z_t|$, the cardinality of random set
$Z_t$. Clearly, the number of multi-edges in $G_t$ is less than
$2mX_t$.

\bl\label{l4} For any $\epsilon>0$, we have
\be\label{3.9}\E(X_t)=O\(t^{\a/2+m(1-\a)+\epsilon}\).\ee\el {\it
Proof}: By the definition of the model, we have
$$\E(X_{t+1}\mid G_{t}) \leq
 X_{t}+(1-\a)\frac{X_{t}}{t}+\a \binom{m}{2}
 \frac{\Delta_{t}}{e_{t}}.$$
Taking expectation and then using (\ref{2.5}) and the fact that
$e_t\geq mt$, we have \bena\E(X_{t+1})&&\hskip-5mm\leq
\(1+\Df{1-\a}{t}\)\E(X_t)+O\(t^{\a/2+m(1-\a)-1}(\log
t)^3\)\nonumber\\[3mm]
&&\hskip-5mm=\(1+\Df{1-\a}{t}\)\E(X_t)+O\(t^{\a/2+m(1-\a)-1+\epsilon}\).\label{3.10}\eena
Using the argument between (\ref{2.18}) and (\ref{2.19}), the Lemma
follows immediately from (\ref{3.10}). \QED

Now, we try to establish the recurrence for $\overline D_k(t)$. Put
$D_{k}(t)=0, 0 \leq k <m$, for all $t\geq 2$. For $k\geq m$, we have
\bena && \overline{D}_{k}(t+1)= \overline{D}_{k}(t)+\a
m\E\(-\Df{kD_k(t)}{2e_t}+\Df{(k-1)D_{k-1}(t)}{2e_t}-O\(\Df{\Delta_t}{e_t}\)\)\nonumber\\[3mm]
&&\hskip 8mm
+(1-\a)(k-1)\E\(-\Df{D_k(t)}{t}+\Df{D_{k-1}(t)}{t}-O\(\Df{X_t}{t}\)\)+\a
\1_{k=m}. \label{3.11}\eena The terms $O\(\Df{\Delta_t}{e_t}\)$ and
$O\(\Df{X_t}{t}\)$ account for the probabilities that we create more
than one degree changes due to new vertex addition and vertex
copying from $Z_t$ respectively.
\vskip 3mm

By Lemma~\ref{l3}, the term $\Dp\E\(\frac{kD_k(t)}{2e_t}\)$ can be
expressed as \bena&&\E\(\left.\frac{kD_k(t)}{2e_t}\right| |e_t-\mu
t|\leq
t^{1/2+\varepsilon_0+m(1-\a)}\)\P\(|e_t-\mu t|\leq t^{1/2+\varepsilon_0+m(1-\a)}\)\nonumber\\[3mm]
&&+\E\(\left.\frac{kD_k(t)}{2e_t}\right| |e_t-\mu t|>
t^{1/2+\varepsilon_0+m(1-\a)}\)\P\(|e_t-\mu t|> t^{1/2+\varepsilon_0+m(1-\a)}\)\nonumber\\[3mm]
&&=\frac{\E\(\left.{kD_k(t)}\right||e_t-\mu t|\leq
t^{1/2+\varepsilon_0+m(1-\a)}\)\P\(|e_t-\mu t|\leq
t^{1/2+\varepsilon_0+m(1-\a)}\)}{2\mu
t}\nonumber\\[2mm]&&\hskip 5mm\times\(1+O\(t^{-1/2+\varepsilon_0+m(1-\a)}\)\)+O(t^{-\varepsilon_0}),\label{3.12}\eena where we used the
fact that $kD_k(t)\leq 2e_t$ to hand the second term. In addition,
we have \bena&&\E\(\left.{kD_k(t)}\right||e_t-\mu t|\leq
t^{1/2+\varepsilon_0+m(1-\a)}\)\P\(|e_t-\mu t|\leq
t^{1/2+\varepsilon_0+m(1-\a)}\)\nonumber\\[2mm]&&=k\overline D_k(t)-\E({kD_k(t)};|e_t-\mu t|>
t^{1/2+\varepsilon_0+m(1-\a)}),\label{3.13}\eena and
\bena&&\E({kD_k(t)};|e_t-\mu t|>
t^{1/2+\varepsilon_0+m(1-\a)})\nonumber\\[3mm]
&&=\E({kD_k(t)};|e_t-\mu t|> t^{1/2+\varepsilon_0+m(1-\a)}, e_t\leq
O(t(\log t)^6))\nonumber\\[2mm]
&&\hskip 5mm+\E({kD_k(t)};|e_t-\mu t|>
t^{1/2+\varepsilon_0+m(1-\a)}, e_t> O(t(\log t)^6))\nonumber\\[3mm]&&\leq
O(t(\log t)^6)\P(|e_t-\mu t|>
t^{1/2+\varepsilon_0+m(1-\a)})\nonumber\\[3mm]
&&\hskip5mm +O(t^2)\P(e_t> O(t(\log t)^6))\nonumber\\[3mm]
&&\leq O(t^{1-\varepsilon_0}(\log
t)^6)+O(t^{-10})=O(t^{1-\varepsilon_0}(\log t)^6).\hskip
25mm\label{3.14}\eena Note that to get (\ref{3.14}), we used the
fact that $kD_k(t)\leq 2e_t$ and the bounds on $e_t$ given in
(\ref{2.22}) and (\ref{2.21}).

Thus, combining (\ref{3.12}), (\ref{3.13}) and (\ref{3.14}),
\bena\E\(\frac{kD_k(t)}{2e_t}\)&&\hskip-4mm=\frac{k\overline
D_k(t)}{2\mu
t}\(1+O\(t^{-1/2+\varepsilon_0+m(1-\a)}\)\)+O(t^{-\varepsilon_0}(\log
t)^6)\nonumber\\[3mm]
&&\hskip-4mm \leq \frac{k\overline D_k(t)}{2\mu t}
+\Df{\E(2e_t)}{2\mu
t}O\(t^{-1/2+\varepsilon_0+m(1-\a)}\)+O(t^{-\varepsilon_0}(\log
t)^6),\nonumber\eena using (\ref{2.9}), we have for $k\geq m$
\bena\E\(\frac{kD_k(t)}{2e_t}\)=\frac{k\overline D_k(t)}{2\mu t}
+O\(t^{-1/2+\varepsilon_0+m(1-\a)}\)+O(t^{-\varepsilon_0}(\log
t)^6).\label{3.15}\eena On the other hand, by inequality (\ref{2.5})
and Lemma~\ref{l4}, for any fixed $\epsilon\in(0,1-\a/2-m(1-\a))$,
we have \be\label{3.16}\E\(\Df{\D_t}{e_t}\),\
\E\(\Df{X_t}{t}\)=O(t^{-1+\a/2+m(1-\a)+\epsilon}).\ee Let
\be\label{3.19}\varepsilon_1=\frac 12\min\left\{\varepsilon_0,1-\a
/2-m(1-\a), 1/2-\varepsilon_0-m(1-\a)\right\}.\ee Now, substitute
(\ref{3.15}) and (\ref{3.16}) into (\ref{3.11}), we get the
recurrence for $\overline{D}_k(t)$ as \bena \overline{D}_{k}(t+1)&&
\hskip -4mm =\overline{D}_{k}(t)-\(\frac
k2-(1-\a)\)\Df{\overline{D}_{k}(t)}{t}+\frac{(k-1)}{2}\Df{\overline{D}_{k-1}(t)}{t}\nonumber\\[3mm]
&&+O(t^{-\varepsilon_1})+\a \1_{k=m},\ \ \forall \ \ k\geq
m.\label{3.17}\eena Note that the hidden constant, denote by $L$, in
term $O(t^{-\varepsilon_1})$ is independent of $k$.

\section{Solving (\ref{3.17}) and The Proof Theorem~\ref{th1}}
\renewcommand{\theequation}{4.\arabic{equation}}
\setcounter{equation}{0}

In recurrence (\ref{3.17}), if we heuristically put $\bar
d_k=\Df{\overline D_k(t)}{t}$ and assume it is a constant, we get
\bena \Dp\frac {(k+2\a)}2\bar d_k&&\hskip-5mm=\frac{(k-1)}{2}\bar
d_{k-1}+O(t^{-\varepsilon_1})+\a\1_{{k=m}}.\nonumber \eena This
leads to the consideration of the following recurrence in $k$:
\be\label{3.2}\left\{\ba{ll}&\hskip-3mm\Df {(k+2\a)}2 d_{k}
=\Dp\frac{(k-1)}{2}d_{k-1}+\a \1_{k=m},\ \ k\geq
m;\\[5mm]
&\hskip-3mmd_k=0,\ \ 0\leq k<m.\ea\right.\ee

The following Lemma shows that (\ref{3.2}) is a good approximation
to (\ref{3.17}).

\bl\label{l5} Suppose that $\{d_k:k\geq 0\}$ be the solution of
(\ref{3.2}), then there exists a constant $M>0$ such that
\be\label{3.3}|\overline D_k(t)-td_k|\leq Mt^{1-\varepsilon_1},\ee
for all $t\geq 1$ and $k\geq 0$, where $\varepsilon_1$ is given in
(\ref{3.19}).\el

{\it Proof}: The recurrence can be solved directly as: $d_k=0$,
$0\leq k<m$; $d_m=\Dp\frac{2\a}{m+2\a}$ and
\be\label{3.18}d_k=\prod_{i=m+1}^k\(1+\Df{1+2\a}{i+2\a}\)d_m,\ \
\forall\ k>m.\ee Obviously, $d_k$ decay as $k^{-(1+2\a)}$,
consequently, for some constant $C$, \be\label{4.1} d_k\leq C/k \hb{
 for all  } k\geq 1.\ee Using (\ref{4.1}) and the degree estimate
given in Lemma~\ref{l1}, the Lemma follows from a standard argument
which can be found in \cite{CFV} (see Lemma 5.1) and \cite{WDLC}
(see Lemma 3.1). \QED

{\it Proof of Theorem~\ref{th1}}: Theorem~\ref{th1} follows
immediately from (\ref{3.3}) and (\ref{3.18}).\QED

\section*{Acknowledgements} The authors thank Prof. Zhao Dong and Prof. Ke Liu for useful advice and discussion.


\begin{thebibliography}{30}

\bibitem{ABJ} R. Albert, A. Barab\'asi and H. Jeong (1999) {\it Diameter of the World Wide
Web}. Nature, {\bf 401}, pp. 103-131.

\bibitem{ACL}
W. Aiello, F. R. K. Chung and L. Lu (2002) {\it Random Evolution in
Massive Graphs} In {\it Handbook on Massive Data Sets}, edited by
James Abello et al., pp. 510-519. Norwood, MA: Kluwer Academic
Publishers

\bibitem{ASBS} L. A. N. Amaral, A. Scala, M. Barth\'el\'emy and H.
E. Stanley (2000) {\it Classes of Small-World Networks,} Proc Natl
Acad Sci U S A. 2000 October 10; {\bf 97}: pp. 11149-11152.

\bibitem{B} B. Bollob\'as (1998) {\it Modern Graph Theory}
Springer-Verlag New York

\bibitem{B2} B. Bollob\'as (2001) {\it Random Graph (second
edition),} Cambridge University Press

\bibitem{BA}A.-L. Barab\'asi and R. Albert (1999) {\it Emergence of Scaling in Random
Networks,} Science {\bf 286}, pp. 509-512


\bibitem{BKEMS}H. R. Bernard, P. D. Killworth, M. J. Evans, C.
McCarty and G. A. Shelley (1988) {\it Studying Social Relations
Cross-Culturally,} Ethnology {\bf 27}, pp. 155-179

\bibitem{BO}
P. G. Buckley and D. Osthus (2004) {\it Popularity Based Random
Graph Model Leading to a Scale-Free Degree Sequence,} Discrete
Mathematics, {\bf 282}, pp. 53-68.


\bibitem{BR} B. Bollob\'as and O. Riordan (2002) {\it Mathematical Results on Scale-Free Random
Graphs.} In {\it Handbook of Graphs and Networks,} pp. 1-34. Berlin:
Wiley-VCH.

\bibitem{BR2} B. Bollob\'as and O. Riordan (2004) {\it The Diameter of a Scale-Free Random Graph,}
Combinatorica {\bf 4}, pp. 5-34.

\bibitem{BR3} B. Bollob\'as and O. Riordan (2003) {\it Robustness and Vulnerability of Scale-Free Random Graph,}
Internet Mathematics {\bf 1}, pp.1-35

\bibitem{BRST} B. Bollob\'as, O. Riordan, J. Spencer and G.
Tusn\'ady (2001) {\it The Degree Sequence of a Scale-Free Random
Graph Process} Random Structure and Algorithms, {\bf 18}, pp.
279-290.

\bibitem{CF}C. Cooper and A. Frieze (2003) {\it A General Model of Undireted Web
Graphs.} Random Structures and Algorithms, {\bf 22}, pp. 311-335.

\bibitem{CFV}
{C. Cooper, A. Frieze and J. Vera (2004)} {\it Random Deletion in a
Scale-Free Random Graph Process.} Internet Mathematics {\bf 1}, pp.
463-483

%

\bibitem{KRRS}R. Kumar, P. Raghavan, S. Rajagopalan, D. Sivakumar, A. Tomkins
and E. Upfal (2000) {\it Stochastic Models for the Web Graph,} In
{\it 41st FOCS}, pp. 57-65.


\bibitem{LLJ}S. Lehmann, B. Lautrup and A. D. Jackson (2003) {\it Citation Networks in High Energy Physics,}
Phys. Rev. E (Statistical, Nonlinear, and Soft Matter Physics), {\bf
68}: 026113

\bibitem{M} M. Mitzenmacher (2003) {\it A Brief History of
Generative Models for Power Law and Lognormal
Distributions,}Internet Mathematics

\bibitem{N} M. E. J. Newman (2003) {\it The Structure and Function of the Complex
Networks,} SIAM Review, {\bf 45}, pp. 167-256.

\bibitem{St}S. H. Strogatz (2001) {\it Exploring Complex Networks,} Nature {\bf 410}, pp.
268-276


\bibitem{WDLC} Xuan-Yuan Wu, Zhao Dong, Ke Liu and Kai-Yuan Cai (2008) {\it On
the Degree Sequence and its Critical Phenomenon of an Evolving
Random Graph Process,} to appear
 \href{http://arxiv.org/PS_cache/arxiv/pdf/0806/0806.4684v1.pdf}
{arXiv:0806.4684v1[math.PR]}

\bibitem{WDLC2} Xuan-Yuan Wu, Zhao Dong, Ke Liu and Kai-Yuan Cai (2008) {\it Phase Transition on The Degree Sequence of a Mixed Random Graph
Process,} to appear
\href{http://arxiv.org/PS_cache/arxiv/pdf/0807/0807.2811v3.pdf}{arXiv:0807.2811v3[math.PR]}

\bibitem{WS}D. J. Watts and S. H. Strogatz (1998) {\it Collective Dynamics of `Small-World' Networks,}
Nature {\bf 393}, pp. 440-442.


\end{thebibliography}
\end{document}